\documentclass[12pt]{article}
\usepackage{amsfonts}
\usepackage{latexsym, amssymb, amsmath, amscd, amsfonts, epsfig, graphicx}
\usepackage{mathrsfs}
\usepackage{color}
\usepackage{ifpdf}
\newtheorem{theorem}{Theorem}[section]

\newtheorem{lemma}[theorem]{Lemma}

\newtheorem{corollary}[theorem]{Corollary}
\newtheorem{conjecture}[theorem]{Conjecture}

 \numberwithin{equation}{section}

\def\qed{\nopagebreak\hfill{\rule{4pt}{7pt}}}
\def\proof{\noindent {\it{Proof.} \hskip 2pt}}

\begin{document}

\begin{center}
{\large\bf The Reverse Ultra Log-Concavity of

the Boros-Moll Polynomials}
\end{center}

\begin{center}
William Y.C. Chen$^{1}$, Cindy C.Y. Gu$^{2}$ \\[6pt]
Center for Combinatorics, LPMC-TJKLC\\
Nankai University, Tianjin 300071, P. R. China

$^{1}${chen@nankai.edu.cn}, $^{2}${guchunyan@cfc.nankai.edu.cn}
\end{center}

\vspace{0.3cm} \noindent{\bf Abstract.} We prove the reverse ultra
log-concavity of the Boros-Moll polynomials. We further establish an
inequality which implies the log-concavity of the sequence
$\{i!d_i(m)\}$ for any $m\geq 2$, where $d_i(m)$ are the
coefficients of the Boros-Moll polynomials $P_m(a)$.  This
inequality also leads to the fact that in the asymptotic sense, the
Boros-Moll sequences are just on the borderline between ultra
log-concavity and  reverse ultra log-concavity.
  We propose two
conjectures on the log-concavity and reverse ultra log-concavity of
the sequence $\{d_{i-1}(m) d_{i+1}(m)/d_i(m)^2\}$ for $m\geq 2$.



\section{Introduction}

This paper is concerned with the reverse ultra log-concavity of the
Boros-Moll polynomials. A sequence $\{a_k\}_{k\geq 0}$ of real
numbers is said to be log-concave if $a_k^2\geq a_{k+1}a_{k-1}$
holds for all $k\geq 1$. A polynomial is said to be log-concave if
the sequence of its coefficients is log-concave, see Brenti
\cite{brenti1989} and Stanley \cite{stanley1989}. Furthermore, a
sequence $ \{a_k\}_{0\leq k\leq n}$ is called ultra log-concave if
$\left\{a_k\big/{n\choose k}\right\}$ is log-concave, see Liggett
\cite{Lig}.  This condition can be restated as
\begin{equation}\label{ultralc}
k(n-k)a_k^2-(n-k+1)(k+1)a_{k-1}a_{k+1}\geq 0.
\end{equation}
It is well known that if a polynomial has only real zeros, then its
coefficients form an ultra log-concave sequence. As noticed by
Liggett \cite{Lig},  if a sequence $\{a_k\}_{0\leq k\leq n}$ is
ultra log-concave, then the sequence $\{k!a_k\}_{0\leq k\leq n}$ is
log-concave.

A sequence is said to be reverse ultra log-concave if it satisfies
the reverse relation of \eqref{ultralc}, that is,
\begin{equation}\label{reverse}
k(n-k)a_k^2-(n-k+1)(k+1)a_{k-1}a_{k+1}\leq 0.
\end{equation}
For example, it is easy to verify that for $n\geq 2$, the Bessel
polynomial \cite{han}
$$y_n(x)=\sum_{k=0}^n\frac{(n+k)!}{2^kk!(n-k)!}x^k$$
is log-concave and  reverse ultra log-concave.

The Boros and Moll polynomials, denoted $P_m(a)$,  arise in the
following evaluation of a quartic integral
$$\int_0^{\infty}\frac{1}{(x^4+2ax^2+1)^{m+1}}dx=\frac{\pi}{2^{m+3/2}(a+1)^{m+1/2}}P_m(a),$$
where \begin{equation} \label{pma}
P_m(a)=2^{-2m}\sum_k2^k{2m-2k\choose m-k}{m+k\choose k}(a+1)^k,
\end{equation}
 see, \cite{Bor2,Bor3,Bor4,Moll}.
Write
$$P_m(a)=\sum_{i=0}^m d_i(m)a^i.$$ The sequence
$\{d_i(m)\}_{0\leq i \leq m}$ is called a Boros-Moll sequence. The
 expression (\ref{pma}) gives the following formula for the
coefficients  $d_i(m)$,
$$d_i(m)=2^{-2m}\sum_{k=0}^m 2^k{2m-2k \choose m-k}{m+k \choose m}{k\choose i}.$$
Clearly, the coefficients $d_i(m)$ are positive. Moll conjectured
that the sequence $\{d_i(m)\}_i$ is  log-concave for $m\geq 2$, that
is, $d_i(m)^2\geq d_{i-1}(m)d_{i+1}(m)$ $(1\leq i \leq m-1)$. This
conjecture has been proved by Kauers and Paule \cite{Pau}.

Despite the log-concavity of $\{d_i(m)\}$,  we find that the inverse
ultra log-concavity holds.

\begin{theorem}\label{new}
For $m\geq 2$ and $1\leq i\leq m-1$, we have
\begin{equation}\label{ulog}
\bigg(\frac{d_{i-1}(m)}{{m\choose
i-1}}\bigg)\cdot\bigg(\frac{d_{i+1}(m)}{{m\choose i+1}}\bigg)>
\bigg(\frac{d_i(m)}{{m\choose i}}\bigg)^2,
\end{equation}
or, equivalently,
\begin{equation}\label{3.6}
\frac{d_i(m)^2}{d_{i-1}(m)d_{i+1}(m)}<\frac{(m-i+1)(i+1)}{(m-i)i}.
\end{equation}
\end{theorem}

On the other hand, it can be shown that the coefficients $d_i(m)$
 satisfy an inequality stronger  than the
log-concavity. To be more specific, we will give a lower bound of
 $d_i(m)^2 /
\left(d_{i-1}(m) d_{i+1}(m)\right)$,  which is very close to the
above upper bound in (\ref{3.6}).

\begin{theorem}\label{c}
For $m\geq 2$ and $1\leq i\leq m-1$, we have
\begin{equation}\label{cin}
\frac{d_i(m)^2}{d_{i-1}(m)d_{i+1}(m)}>
\frac{(m-i+1)(i+1)(m+i)}{(m-i)i(m+i+1)}.
\end{equation}
\end{theorem}

This paper is organized as follows. We  establish an upper bound of
$d_i(m+1)/d_i(m)$ in Section \ref{sec3}, which leads to the reverse
ultra log-concavity of $\{d_i(m)\}$. In Section \ref{sec5} we give
the proof of Theorem \ref{c}. We conclude this paper with two
conjectures concerning  the log-concavity and the reverse ultra
log-concavity of the sequence  $\{d_{i-1}(m)d_{i+1}(m)/d_i^2(m)\}$
for $m\geq 2$.

\section{An Upper Bound for $d_i(m+1)/d_i(m)$}\label{sec3}

In this section, we  establish an upper bound for the ratio
$d_i(m+1)/d_i(m)$ that will lead to the reverse ultra log-concavity
of the sequence of $\{d_i(m)\}$. For $m\geq 1$ and $0\leq i \leq m$,
set
\begin{equation}\label{3.2}
T(m,i)=\frac{4m^2+7m+3+i\sqrt{4m+4i^2+1}-2i^2}{2(m-i+1)(m+1)}.
\end{equation}

\begin{theorem}\label{lem3.1}
For all $m\geq 2,\,1\leq i\leq m-1$, we have
\begin{equation}\label{3.1}
{d_i(m+1)\over d_i(m)}< T(m,i),
\end{equation}
and for  $m\geq 1$, we have
\begin{equation}\label{d0}
{d_0(m+1)\over d_0(m)}=T(m,0),\,\quad {d_m(m+1)\over d_m(m)}=
T(m,m).
\end{equation}
\end{theorem}

The following lemma will be needed in the proof of Theorem
\ref{lem3.1}.

\begin{lemma}\label{tf}
For $m\geq 2$ and $1\leq i\leq m-1$,
\begin{equation}\label{3.4}
T(m,i)< F(m,i),
\end{equation}
where
$$F(m,i)=\frac{(m+i+1)(4m+3)(4m+5)}{2(m+1)(4m^2-2i^2+9m+5-i\sqrt{4m+4i^2+5})}.$$
\end{lemma}

\proof Let $A=\sqrt{4m+4i^2+1}$ and $B=\sqrt{4m+4i^2+5}$. It is easy
to check that
\begin{equation}\label{fsubt}
F(m,i)-T(m,i) =\frac{i(X-Y)}{2(m+1)(m-i+1)(4m^2+9m+5-2i^2-iB)},
\end{equation}
where
\begin{align*} X= & (i-4i^3)+iAB\\[3pt]
 Y= & (5+4m^2+9m-2i^2)A-(3+4m^2+7m-2i^2)B.
\end{align*}
Since $(4m^2+9m+5-2i^2)^2-(iB)^2=(4m+5)^2(m+i+1)(m-i+1)> 0,$ it
remains to show that the numerator of \eqref{fsubt} is also
positive. We claim that $X>0$ and $X^2>Y^2$.

 Since $m> i$, we have $A>2i+1$
and $B>2i+1$. Moreover, since $i\geq 1$, we find that
$$X=(i-4i^3)+iAB\geq i-4i^3+i(2i+1)^2=4i^2+2i> 0.$$
It is routine to check $X^2-Y^2=G(m,i)-H(m,i)$, where
\begin{eqnarray*}
G(m,i)& = &
(32m^4-32m^2i^2+128m^3-64mi^2+190m^2-30i^2+124m+30)AB,\\[3pt]
H(m,i)& =&128m^5+608m^4+1128m^3+1014m^2+436m+128m^4i^2+384m^3i^2\\[3pt]
&&  +\,408m^2i^2 -128m^2i^4+200mi^2-256mi^4-120i^4+50i^2+70.
\end{eqnarray*}
Since $i < m$, it is easily seen that $G(m,i)>0$ and $H(m,i)>0$. To
prove $G(m,i)>H(m,i)$, it suffices to show that $G(m,i)^2
>H(m,i)^2$. In fact, for $1\leq i\leq m-1$,
$$G(m,i)^2-H(m,i)^2=16(4m+5)^2(16mi^2+12i^2-1)(m+i+1)^2(m-i+1)^2> 0.$$
This yields $X^2>Y^2$. Since $X>0$, we see that $X>Y$, and hence
\eqref{3.4} holds for $1\leq i\leq m-1$.\qed

\vspace {6mm}

 \noindent{\emph{Proof of Theorem
\ref{lem3.1}.}} It is easy to check (\ref{d0}). To prove
\eqref{3.1}, we  proceed by induction on $m$. For $m=2$ and $i=1$,
we have $d_{1}(3) / d_1(2)= 43/15< T(2,1)=(31+\sqrt{13})/12$. We now
assume that \eqref{3.1}  is true for $m$, that is,
\begin{equation}\label{induc1}
d_i(m+1)<T(m,i)d_i(m),\quad  1\leq i\leq m-1.
\end{equation}
It will be shown  that
\begin{equation}\label{induc2}
d_i(m+2)<T(m+1,i)d_i(m+1),\quad 1\leq i\leq m-1.
\end{equation}

Using the recurrence \eqref{2.3}, we may write \eqref{induc2} in the
following form
\begin{align}
 \frac{-4i^2+8m^2+24m+19}{2(m-i+2)(m+2)}d_i(m+1)
&-\frac{(m+i+1)(4m+3)(4m+5)}{4(m+1)(m+2)(m-i+2)}d_i(m) \nonumber
\\[5pt]&<
T(m+1,i)d_i(m+1).\label{2.7}
\end{align}
Since $m> i$, we have $4m+4i^2+5<12m +4m^2+9$. It follows that
\begin{align*}
R(m,i)=&\,\frac{-4i^2+8m^2+24m+19}{2(m-i+2)(m+2)}-T(m+1,i)\\[5pt]
=&\,\frac{4m^2+9m+5-2i^2-i\sqrt{4m+4i^2+5}}{2(m-i+2)(m+2)}\\[5pt]
\geq\, &\frac{4m^2+9m+5-2i^2-i(2m+3)}{2(m-i+2)(m+2)}>0.
\end{align*}
Therefore, \eqref{2.7} is equivalent to the following inequality
\begin{equation}\label{dimr}
{d_i(m+1)\over d_i(m)}< F(m,i),
\end{equation}
which is a consequence of \eqref{induc1} and Lemma \ref{tf}.

It remains to consider the case $i=m$. We aim to show that
\begin{equation} \label{dt}
\frac{d_m(m+2)}{d_m(m+1)} < T(m+1,m).
\end{equation}
By easy computation, we find that
\begin{align*}
\frac{d_m(m+2)}{d_m(m+1)}&=\frac{(m+1)(4m^2+18m+21)}{2(2m+3)(m+2)},\\[3pt]
T(m+1,m)&=\frac{2m^2+15m+14+m\sqrt{4m^2+4m+5}}{4(m+2)}.
\end{align*}
Thus  (\ref{dt}) can be rewritten as
\begin{equation} \label{2m4m}
(2m^2+3m)\sqrt{4m^2+4m+5} > 4m^3+8m^2+5m.
\end{equation}
Denote by $U$ and $V$ the left hand side and the right hand side of
(\ref{2m4m}), respectively. Then, $U^2-V^2= 4m^2(4m+5)>0$, and so
(\ref{dt}) is verified. This  completes the proof. \qed

\section{The Reverse Ultra Log-concavity}\label{sec4}

In this section, we give the proof of Theorem \ref{new}. Our
approach can be described as follows.  Let $f(x)=ax^2+bx+c$ be a
quadratic function with $a>0$. Suppose that the equation $f(x)=0$
has two distinct real zeros $x_1$ and $x_2$, where $x_1<  x_2$. Then
$f(x)> 0$ if $x>  x_2$ or $x<  x_1$ and $f(x)< 0$ if $x_1<  x< x_2$.
The key step is to transform the inequality \eqref{3.6}, that is,
\[
\frac{d_i(m)^2}{d_{i-1}(m)d_{i+1}(m)}<\frac{(m-i+1)(i+1)}{(m-i)i},
\]
into a quadratic inequality in the ratio $d_i(m+1)/d_i(m)$.

We will need the following recurrence relations for the coefficients
$d_i(m)$. For $m\geq 1$ and $0\leq i\leq m+1$,
\begin{align}
  2(m+1)d_i(m+1)&=2(m+i)d_{i-1}(m)+(4m+2i+3)d_i(m),
\label{2.1}\\[5pt]
  2(m+1)(m+1-i)d_i(m+1)&=(4m-2i+3)(m+i+1)d_i(m) \nonumber\\[5pt]
  &\quad -2i(i+1)d_{i+1}(m),
\label{2.2}\\[5pt]
4(m+2-i)(m+1)(m+2)d_i(m+2)
&=2(m+1)(-4i^2+8m^2+24m+19)d_i(m+1)\nonumber\\[5pt]
&\quad-(m+i+1)(4m+3)(4m+5)d_i(m). \label{2.3}
\end{align}
These recurrences are derived by Kauers and Paule \cite{Pau}. The
relation \eqref{2.3} is also derived independently by Moll
\cite{Moll 2}. Based on these recurrence relations, Kauers and Paule
\cite{Pau} derived the following lower bound of $d_i(m+1)/d_i(m)$
 in their proof of the log-concavity of Boros-Moll
polynomials
\begin{equation}\label{2.5}
{d_i(m+1)\over d_i(m)}\geq Q(m,i),\quad 0\leq i\leq m,
\end{equation}
where
\begin{equation}\label{qmi}
Q(m,i)=\frac{4m^2+7m+i+3}{2(m+1-i)(m+1)}.
\end{equation}
Note that Chen and Xia \cite{chen} have shown that the above
inequality  \eqref{2.5} becomes strict for $m\geq 2$ and $1\leq i
\leq m-1$, that is,
\begin{equation}\label{q}
{d_i(m+1)\over d_i(m)}> Q(m,i).
\end{equation}

Now we are ready to prove the reverse ultra log-concavity of
$\{d_i(m)\}$.

 \noindent{\emph{Proof of Theorem \ref{new}.}}
Applying
 \eqref{2.1} and \eqref{2.2}, we may reformulate  \eqref{3.6} in the
 following form
\begin{align}
&\ 4(m-i+1)^2(m+1)^2\bigg(\frac{d_i(m+1)}{d_i(m)}\bigg)^2 \nonumber\\[3pt]
&-4(m-i+1)(m+1)(4m^2-2i^2+7m+3)\bigg(\frac{d_i(m+1)}{d_i(m)}\bigg)\nonumber\\[3pt]
&-(32mi^2-56m^3-73m^2-42m+13i^2-9-16m^4+16i^2m^2)< 0. \label{3.8}
\end{align}
For $1\leq i \leq m-1$, the discriminant of the above quadratic
function in $d_i(m+1)/d_i(m)$ equals
$$\bigtriangleup=16i^2(m+1)^2(4i^2+4m+1)(m-i+1)^2> 0.$$
We see that  the quadratic function on the left hand side of
\eqref{3.8} has two real roots
\begin{align*}
x_1&=\frac{4m^2-2i^2+7m+3-i\sqrt{4m+4i^2+1}}{2(m-i+1)(m+1)}, \\[3pt]
x_2&=\frac{4m^2-2i^2+7m+3+i\sqrt{4m+4i^2+1}}{2(m-i+1)(m+1)}.
\end{align*}
 Clearly, $Q(m,i)>x_1$. In view of \eqref{2.5}, we deduce that
  $d_i(m+1)/d_i(m)\geq Q(m,i)>x_1$. Observe that $x_2$ coincides
  with the upper bound $T(m,i)$ in Theorem \ref{lem3.1}.
  Thus we have $d_i(m+1)/d_i(m)<x_2$.
  So we have shown that  for $1\leq i\leq m-1$,
$$x_1<{d_i(m+1)\over d_i(m)}<x_2, $$
which implies \eqref{3.8}. This completes the proof of Theorem
\ref{new}. \qed

\section{A Lower Bound for $d_i(m)^2/(d_{i-1}(m)d_{i+1}(m))$}\label{sec5}

In this section, we give the proof of Theorem \ref{c} on a lower
bound of  $d_i(m)^2/(d_{i-1}(m)d_{i+1}(m))$. As will be seen, the
lower bound for $d_i(m)^2/(d_{i-1}(m)d_{i+1}(m))$ is very close to
the upper bound (\ref{3.6}) for the reverse ultra log-concavity. So
in the asymptotic sense, we may say that the Boros-Moll polynomials
are just on the borderline between ultra log-concavity and reverse
ultra log-concavity. We conclude this paper with two conjectures.

\noindent{\emph{Proof of Theorem \ref{c}.}} Utilizing the recurrence
relations (\ref{2.1}) and (\ref{2.2}), the inequality \eqref{cin}
can be restated as
\begin{align*}
&\
4(m+1)^2(m-i+1)^2\left(\frac{d_i(m+1)}{d_i(m)}\right)^2\\[3pt]
&-4(m-i+1)(m+1)(4m^2+7m-2i^2+3)\frac{d_i(m+1)}{d_i(m)}\\[3pt]
&+(4m^2+7m+3)(-4i+3+4m)(m+i+1)> 0.
\end{align*}
For $1\leq i\leq m-1$, the discriminant of the above quadratic
function in $d_i(m+1)/d_i(m)$ equals
\begin{equation}
\delta=16i^2(2i+1)^2(m+1)^2(m-i+1)^2>0.
\end{equation}
Hence the above quadratic function has two real roots,
\begin{align*}
x_1&=\frac{4m^2+7m-4i^2-i+3}{2(m+1)(m-i+1)},\\[3pt]
x_2&=\frac{4m^2+7m+i+3}{2(m+1)(m-i+1)}.
\end{align*}
As $x_2=Q(m,i)$, it follows from  \eqref{q} that
$d_i(m+1)/d_i(m)>x_2$. So we arrive at \eqref{cin}. This completes
the proof. \qed

Notice that for $1\leq i\leq m-1$,
$$\frac{(m-i+1)(i+1)(m+i)}{(m-i)i(m+i+1)}>\frac{i+1}{i}.$$
 As a consequence of Theorem \ref{c} , we obtain the log-concavity
of the sequence $\{i! d_i(m)\}$.

\begin{corollary}\label{dll}
For $m\geq 2$ and $1\leq i\leq m-1$,
\begin{equation}\label{log2}
\frac{d_i^2(m)}{d_{i-1}(m)d_{i+1}(m)}> \frac{i+1}{i},
\end{equation}
or equivalently, the sequences $\{i!d_i(m)\}$ is log-concave.
\end{corollary}

\begin{corollary}\label{g} For $1\leq i\leq m-1$, let
$$  c_i(m)=\frac{d_i^2(m)}{d_{i-1}(m)d_{i+1}(m)} \quad
\mbox{and}\quad
u_i(m)=\left(1+\frac{1}{i}\right)\left(1+\frac{1}{m-i}\right).$$
Then for any $i\geq 1$,
\begin{equation} \label{limit-1}
\lim_{m\rightarrow\infty}\frac{c_i(m)}{u_i(m)}=1.
\end{equation}
\end{corollary}

\proof By Theorems \ref{new} and \ref{c}, we find that
$$\frac{m+i}{m+i+1}< \frac{c_i(m)}{u_i(m)}<1,$$
which implies (\ref{limit-1}). \qed

We remark that even when  $m$ is small, $c_i(m)$ is quite close to
$u_i(m)$ for any $1\leq i\leq  m-1$. Numerical evidence indicates
that $c_i(m)/u_i(m)$ is increasing for given $m$. For example, when
$m=8$, the values of $c_i(m)/u_i(m) $ for $1\leq i\leq 7$ are given
below
\[ 0.956593,\quad 0.969751, \quad
0.978293, \quad  0.983956, \quad 0.987811, \quad 0.990507, \quad
0.992445.
\]

We propose the following two conjectures on the log-concavity and
reverse ultra log-concavity of the sequence $\{  d_{i+1}(m)
d_{i-1}(m) /d_i(m)^2\}$.
\begin{conjecture}  \label{conj-2}
For $m\geq 2$, the sequence $\{  d_{i+1}(m) d_{i-1}(m)
/d_i(m)^2\}_{2\leq i \leq m-2}$ is log-concave.
 \end{conjecture}

\begin{conjecture}\label{conj-1}
For $m\geq 2$, the sequence $\{  d_{i+1}(m) d_{i-1}(m)
/d_i(m)^2\}_{2\leq i\leq m-2}$ is reverse ultra log-concave.
\end{conjecture}

\vspace{.2cm} \noindent{\bf Acknowledgments.} The authors would like
to thank the referee for the  valuable comments leading to an
improvement of an earlier version. This work was supported by the
973 Project, the PCSIRT Project of the Ministry of Education, the
Ministry of Science and Technology, and the National Science
Foundation of China.

\end{document}